# Infinite set of non linear Equations for the Li- Keiper Coefficients: a possible new upper and lower bound

Danilo Merlini[1, 2, 3], Massimo Sala[2, 3], Nicoletta Sala[2, 3]
Cerfim, Switzerland
ISSI, Switzerland
IRFIA, Switzerland
(Research Institute in Arithmetic Physics)

Abstract

Starting with an infinite set of non linear Equations for the Li-Keiper coefficients, we first specify a lower bound emerging from the infinite set and give a characterization of it. Then, we propose a possible new upper and lower bound for the coefficients in few of the partitions occurring in the cluster functions furnishing in a nonlinear way the coefficients.
A numerical experiment up to n=15 confirms the proposed bounds and an experiment, i.e. the counting of the zeros in the binary representation of an integer for a constant related to the Glaisher-Kinkelin constant is also given up to n=32.

Plan of the note

1. Introduction
2. An infinite set of equations and a numerical check.
3. Upper and lower bound for the low Li-Keiper coefficients and special solutions: Riemann Wave Background .
4. Shift of one unity in Xi(s) and the Glaisher-Kinkelin constant.
5. Concluding remark and Riemann Hypothesis.

## 1. Introduction

In this work we analyse an alternative approach to our recent work concerning a possible proof of the Riemann Hypothesis using Theorems of Statistical Mechanics for a model in dimension 1, with long range interaction related to a polynomial truncation of the Xi function in the variable $z = 1 - 1/s$ (the zeros in z are on the unit circle).

The alternative approach consists in the analysis of an infinite set of Equations (called equilibrium equation) for the partial partition functions furnishing the Li- Keiper coefficient.

## 2. An infinite set of equations and a numerical check

Recently, we have obtained an infinite set of nonlinear Equations for the Li-Keiper coefficients, by means of their definition. They are different from other known Equations and involve here the partial partition function $\varphi_n$ at the left hand side and the successive derivatives of the function at the right border of the critical strip: these could be called the "equilibrium equations" for the Li-Keiper coefficients. They read [1]

$$\sum_{k=0}^{n-1} (-1)^k \cdot \binom{n-1}{k} \cdot \varphi_{n-k} = \left( \frac{\xi^n(1)}{\Gamma(n+1)} \right)$$

(1)

where $\varphi_n$ is the coefficient of $z^n$ in the expansion of the partition function, i.e. the Xi function given by

$2 \cdot \xi(z) = \exp(\sum (\lambda_n/n) \cdot z^n) = \sum \varphi_n \cdot z^n =$

$= \varphi_0 + \varphi_1 \cdot z + \varphi_2 \cdot z^2 + \varphi_3 \cdot z^3 + ... =$

$= 1 + \lambda_1 \cdot z + (\frac{1}{2}) \cdot (\lambda_1^2 + \lambda_2) \cdot z^2 + (1/3) \cdot (\lambda_3 + (3/2) \cdot \lambda_1 \cdot \lambda_2 + (1/2) \cdot \lambda_1^3) \cdot z^3 + ...$

and
$\xi^n(1) = (d^n/ds^n (\xi(s))) |_{s=1}$.

The first few six Equations of the set are:

$$1 \cdot \varphi_1 = (\xi^1/\xi)|_{s=1}/1!$$
$$1 \cdot \varphi_2 - 1 \cdot \varphi_1 = (\xi^2/\xi)|_{s=1}/2!$$
$$1 \cdot \varphi_3 - 2 \cdot \varphi_2 + 1 \cdot \varphi_1 = (\xi^3/\xi)|_{s=1}/3!$$
$$1 \cdot \varphi_4 - 3 \cdot \varphi_3 + 3 \cdot \varphi_2 - 1 \cdot \varphi_1 = (\xi^4/\xi)|_{s=1}/4!$$
$$1 \cdot \varphi_5 - 4 \cdot \varphi_4 + 6 \cdot \varphi_3 - 4 \cdot \varphi_2 + 1 \cdot \varphi_1 = (\xi^5/\xi)|_{s=1}/5!$$
$$1 \cdot \varphi_6 - 5 \cdot \varphi_5 + 10 \cdot \varphi_4 - 10 \cdot \varphi_3 + 5 \cdot \varphi_2 - 1 \cdot \varphi_1 = (\xi^6/\xi)|_{s=1}/6!$$

(2)

As a check, we compute the left hand side of the Equations in calculating the $\varphi$'s using the values (up to some decimals) of the lambda's given in [2] and for the right hand side using the Riemann series at $s=1/2$ or simply the expansion with Maple of the Taylor series of Xi(s) at $s=1$ (Table 1).

| Left | Right |
|---|---|
| 0.0230957089 | 0.0230957089 |
| 0.0233438645 | 0.0233438645 |
| 0.0004979838 | 0.0004979846 |
| 0.0002531817 | 0.0002531898 |
| 0.0000050504 | 0.0000050493 |
| 0.0000017212 | 0.0000172087 |

**Table 1**

### 3. Upper and lower bound for the low Li-Keiper coefficients and special solutions: Riemann Wave Background

In the above Equations we set $\varphi_k = k \cdot \varphi_1 = k \cdot \lambda_1$ for every positive integer k; thus, for the left hand side of Eq.(1) we have:

$$\sum_{k=0}^{n-1}(-1)^k \cdot \binom{n-1}{k} \cdot (n-k) \cdot \lambda_1 = 0$$

(3)

for every n.
Thus:
$$0 '=' (\xi^n(1)/\xi(1))/n! > 0 \quad (4)$$
In fact relation (4) is to understand as follows:

From $\varphi_2 - \varphi_1 = (\xi^2/\xi)|_{s=1}/2!$ , we have that $\varphi_2 = 2 \cdot \varphi_1 + \delta_2$
with $\delta_2 = (\xi^2/\xi)|_{s=1}/2! - \varphi_1$, with $\delta_2 > 0$ [1, 3].
Thus
$$\varphi_2 = 2 \cdot \varphi_1 + \delta_2 \text{ with } \delta_2 > 0.$$

Next from:
$\varphi_3 - 2\varphi_2 + \varphi_1 = (\xi^3/\xi)|_{s=1}/3! = \delta_3 > 0$
we have:

$\varphi_3 = 2\varphi_2 - \varphi_1 + \delta_3 = 2 \cdot (2 \cdot \varphi_1 + \delta_2) - \varphi_1 + \delta_3 = 3 \cdot \varphi_1 + 2 \cdot \delta_2 + \delta_3$

Thus: $\varphi_3 = 3 \cdot \varphi_1 + 2 \cdot \delta_2 + \delta_3$ with $2 \cdot \delta_2 + \delta_3 > 0$.

From
$\varphi_4 = 3 \cdot \varphi_3 - 3 \cdot \varphi_2 + \varphi_1 + (\xi^4/\xi)|_{s=1}/4! = 3 \cdot \varphi_3 - 3 \cdot \varphi_2 + \varphi_1 + \delta_4$
with $\delta_4 = (\xi^4/\xi)|_{s=1}/4! > 0$,
$\varphi_4 = 3 \cdot (3 \cdot \varphi_1 + 2 \cdot \delta_2 + \delta_3) - 3 \cdot (2 \cdot \varphi_1 + \delta_2) + \varphi_1 + \delta_4$.

Thus:
$\varphi_4 = 4 \cdot \varphi_1 + 3 \cdot \delta_2 + 3 \cdot \delta_3 + \delta_4$ with $3 \cdot \delta_2 + 3 \cdot \delta_3 + \delta_4 > 0$.

From $\varphi_5 = 4 \cdot \varphi_4 - 6 \cdot \varphi_3 + 4 \cdot \varphi_2 - \varphi_1 + (\xi^5/\xi)|_{s=1}/5!$

In the same way, with $\delta_5 = (\xi^5/\xi)|_{s=1}/5! > 0$, we obtain
$\varphi_5 = 5 \cdot \varphi_1 + 4 \cdot \delta_2 + 6 \cdot \delta_3 + 4 \cdot \delta_4 + \delta_5 = 5 \cdot \varphi_1 + \varepsilon_5$
with $\varepsilon_5 > 0$.
Finally, from $\varphi_6 - 5 \cdot \varphi_5 + 10 \cdot \varphi_4 - 10 \cdot \varphi_3 + 5 \cdot \varphi_2 - 1 \cdot \varphi_1 = (\xi^6/\xi)|_{s=1}/6!$
We obtain:
$\varphi_6 = 6 \cdot \varphi_1 + 5 \cdot \delta_2 + 10 \cdot \delta_3 + 10 \cdot \delta_4 + 5 \cdot \delta_5 + \delta_6 = 6 \cdot \varphi_1 + \varepsilon_6$

with $\varepsilon_6 = 5 \cdot \delta_2 + 10 \cdot \delta_3 + 10 \cdot \delta_4 + 5 \cdot \delta_5 + \delta_6 > 0$.

Then:
$$\varphi_n = n \cdot \varphi_1 + \varepsilon_n = n \cdot \lambda_1 + \varepsilon_n, \text{ with } \varepsilon_n > 0 \qquad (5)$$
(thanks to the Riemann Series) and $\varepsilon_n$ is increasing with n, and
$\lambda_1 = 1 + \gamma/2 - (½) \cdot \log(4\pi) = 0.0230957...$
where $\gamma$ is the Euler-Mascheroni constant, Wolf [4].

Moreover we obtain:
$$\varphi_n = n \cdot \varphi_1 + \sum_{k=1}^{n-1} \binom{n-1}{k} \cdot \delta_{k+1}$$
(6)

Our rigorous positive lower bound reads:
$$\varphi_n > n \cdot \lambda_1 \text{ for all } n > 1. \qquad (7)$$

Concerning the fluctuation of the $\varphi_s$ of Eq.(5) there is a binomial structure given by:
$$\varepsilon_n = \sum_{k=1}^{n-1} \binom{n-1}{k} \cdot \delta_{k+1}$$
(8)

$$\varphi_n = n \cdot \varphi_1 + \sum_{k=0}^{n-1} \binom{n-1}{k} \cdot \delta_{k+1}$$
(9)

Explicitly:
$$\varphi_n = n \cdot \varphi_1 + (n-1) \cdot \delta_2 + \sum_{k=2}^{n-1} \binom{n-1}{k} \cdot \frac{\frac{\xi^{k+1}}{\xi}\big|_{s=1}}{(k+1)!}$$
(10)

or, since from above we defined $\delta_2 = (\xi^2/\xi)|_{s=1}/2! - \varphi_1$, where $\varphi_1 = (\xi^1/\xi)|_{s=1} = \lambda_1$, we have:

$$\varphi_n = \frac{\xi^1}{\xi}\bigg|_{s=1} + (n-1) \cdot \frac{\frac{\xi^2}{\xi}\big|_{s=1}}{2!} + \sum_{k=2}^{n-1} \binom{n-1}{k} \cdot \frac{\frac{\xi^{k+1}}{\xi}\big|_{s=1}}{(k+1)!}$$

$$\varphi_n = \sum_{k=0}^{n-1} \binom{n-1}{k} \cdot \frac{\frac{\xi^{k+1}}{\xi}\big|_{s=1}}{(k+1)!}$$
(11)

Remark

We notice here that the $\varphi_n$ are related to the complete Bell Polynomials or to the Faà di Bruno's formula [5] in the form of determinants. We recall here the formula for $\varphi_3$.

The known general formula is:

$$exp\left(\sum_{}^{\infty} a_j \cdot \frac{z^j}{j!}\right) = \sum_{k=0}^{\infty} B_k(a_1..a_2..a_k) \cdot \frac{z^k}{k!}$$

where $B_k$ is the complete k-Bell Polynomial and where the Li-Keiper coefficients are given by $\lambda_j = a_j/(j-1)!$, the corresponding determinant for n=3 reads:

$$B_3(a_1, a_2, a_3) = Det \begin{pmatrix} a_1 & 2 \cdot a_2 & a_3 \\ -1 & a_1 & a_2 \\ 0 & -1 & a_1 \end{pmatrix} = a_1^3 + a_1 \cdot a_2 + 2 \cdot a_1 \cdot a_2 + a_3 =$$
$$= \lambda_1^3 + 3 \cdot \lambda_1 \cdot \lambda_2 + 2 \cdot \lambda_3$$

since $a_3 = 2 \cdot \lambda_3$, then:

$(1/3!) \cdot (\lambda_1^3 + 3 \cdot \lambda_1 \cdot \lambda_2 + 2 \cdot \lambda_3) = (1/3) \cdot (\lambda_3 + (3/2) \cdot \lambda_1 \cdot \lambda_2 + \lambda_1^3/2) = \varphi_3$.

We now continue and first check the Formula above for n=5 for which we know, inserting the first 5 lambda 's values taken from Ref [2], that $\varphi_5 = 0.12047684...$
From Eq.(9) or Eq.10 with n=5, we obtain $\varphi_5 = 0.12047684..$ in full agreement with the above number.
For an analysis of subsequent Expressions we analyse numerically the behavior of $\varphi_n$ as a function of n, calculated with Eq.(11) at s=1 that is a Formula with the absence of non trivial zeros; the results up to n= 10 are given in the Table 2 and in the plot (notice: $\xi(s) = \sum_{n=0} a_n \cdot (s-1)^n$).

$$\varphi_n = \sum_{k=0}^{n-1} \binom{n-1}{k} \cdot \frac{\left.\frac{\xi^{k+1}}{\xi}\right|_{s=1}}{(k+1)!} = \sum_{k=0}^{n-1} \binom{n-1}{k} \cdot 2 \cdot a_{k+1}$$

(12)

| n | $\varphi_n$ |
|---|---|
| 1 | 0.0230957088 |
| 2 | 0.0464395736 |
| 3 | 0.0702814232 |
| 4 | 0.0948744395 |
| 5 | 0.1204768540 |
| 6 | 0.1473536683 |
| 7 | 0.1757784078 |
| 8 | 0.2060349171 |
| 9 | 0.2384192110 |
| 10 | 0.2732413926 |

**Table 2**

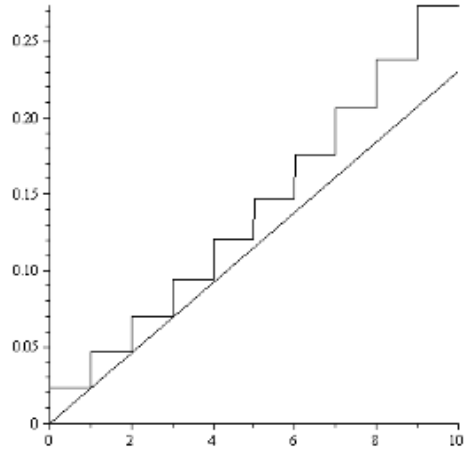

Figure 1. Plot of $\varphi_n(n)$ and the lower bound $n \cdot \varphi_1 = n \cdot \lambda_1$

A characterization of the above special "solution" $\varphi_n = n \cdot \varphi_1$, furnishes, after summation over n in Eq.(2), the function [1]

$$\exp(\log(1+ \varphi_1 \cdot z/(1-z)^2)) = \exp(\log(1+ \varphi_1 \cdot K(z))) \qquad (13)$$

where K(z) is the Koebe function of argument z. Moreover the expansion of the above log function i.e. the coefficient of $z^n$ in the expansion of $\log(1+\lambda_1 \cdot z/(1-z)^2)$ (up to the factor 1/n) is our lower bounds for each n on the Li-Keiper coefficients and such a sequence constitutes our Riemann Wave Background [1]. The series is given explicitly by:

$$\exp(\log(1+ \lambda_1 \cdot z/(1-z)^2)) = \lambda_1 \cdot z + (1/2) \cdot (4 \cdot \lambda_1 - \lambda_1^2) \cdot z^2 + (1/3) \cdot (9 \cdot \lambda_1 - 6 \cdot \lambda_1^2 + (1/2) \cdot \lambda_1^3) \cdot z^3 +$$

$$+ (1/4) \cdot (16 \cdot \lambda_1 - 20 \cdot \lambda_1^2 + 8 \cdot \lambda_1^3 - \lambda \cdot 14) \cdot z^4 + \ldots + \ldots \qquad (14)$$

We may also consider the Formula given by Eq.(10) to obtain the general formula for the Li-Keiper coefficients.
We have:

$$\varphi_1 = (\xi^1/\xi)|_{s=1} = \lambda_1 = 1 \cdot \varphi_1$$

From $\varphi_2 = (\frac{1}{2}) \cdot (\lambda_2 + \lambda_1^2) \rightarrow \lambda_2 = 2 \cdot \varphi_2 - 1 \cdot \varphi_1 \cdot \varphi_1$

From $\varphi_3 = (1/3) \cdot (\lambda_3 + (3/2) \cdot \lambda_1 \cdot \Lambda_2 + \lambda_1^3/2) \rightarrow$
$\lambda_3 = 3 \cdot \varphi_3 - (3/2) \cdot \varphi_1 \cdot (2 \cdot \varphi_2 - 1 \cdot \varphi_1 \cdot \varphi_1) - \varphi_1^3/2 = 3 \cdot \varphi_3 - 3 \cdot \varphi_1 \cdot \varphi_2 + 1 \cdot \varphi_1 \cdot \varphi_1 \cdot \varphi_1.$

From $\varphi_4 = (\frac{1}{4}) \cdot (\lambda_4 + (4/3) \cdot \lambda_1 \cdot \lambda_3 + (\frac{1}{2}) \cdot \lambda_2^2 + \lambda_1^2 \cdot \lambda_2 + \lambda_1^4/6) \rightarrow$
$\lambda_4 = 4 \cdot \varphi_4 - (4 \cdot \varphi_1 \cdot \varphi_3 + 2 \cdot \varphi_2 \cdot \varphi_2) + 4 \cdot \varphi_1 \cdot \varphi_1 \cdot \varphi_2 - 1 \cdot \varphi_1 \cdot \varphi_1 \cdot \varphi_1 \cdot \varphi_1.$

From $\varphi_5 = (1/5) \cdot (\lambda_5 + (5/4) \cdot \lambda_1 \cdot \lambda_4 + (5/6) \cdot \lambda_2 \cdot \lambda_3 + (5/6) \cdot \lambda_1^2 \cdot \lambda_3 +$
$+ (5/8) + \lambda_1 \cdot \lambda_2^2 + (5/12) \cdot \lambda_1^3 \cdot \lambda_2 + (5/120) \cdot \lambda_1^5$

we have:
$\lambda_5 = 5 \cdot \varphi_5 - (5 \cdot \varphi_1 \cdot \varphi_4 + 5 \cdot \varphi_2 \cdot \varphi_3) + (5 \cdot \varphi_1 \cdot \varphi_1 \cdot \varphi^3 + 5 \cdot \varphi_1 \cdot \varphi_2 \cdot \varphi_2) - 5 \cdot \varphi_1 \cdot \varphi_1 \cdot \varphi_1 \cdot \varphi_2 +$
$+ 1 \cdot \varphi 1 \cdot \varphi 1 \cdot \varphi 1 \cdot \varphi 1 \cdot \varphi 1$

i.e. here (5, -10, +10, -5, +1).

Moreover, from φ$_6$= ... we have:

λ$_6$ = 6·φ$_6$ - (6·φ$_1$·φ$_5$+6· φ$_2$·φ$_4$ +3· φ$_3$·φ$_3$ )+(6·φ$_1$·φ$_1$·φ$_4$+12·φ$_1$·φ$_2$·φ$_3$ +
 +2·φ$_2$·φ$_2$·φ$_2$ )¬ (9·φ$_1$·φ$_1$·φ$_2$·φ$_2$ + 6·φ$_1$·φ$_1$·φ$_1$·φ$_3$)+ (6·φ$_1$·φ$_1$·φ$_1$·φ$_1$·φ$_2$) ¬
 ¬1·φ$_1$·φ$_1$·φ$_1$·φ$_1$·φ$_1$·φ**$_1$** ,

i.e. here ( 6, -15,+ 20,-15,+6, -1).

We then have for each n:

$$\lambda_n = \sum_{k=1}^{n}(-1)^{(k-1)} \cdot \left( \sum_{\sum_{i=1}^{k} p_i = n} a_k(p_1, p_2, \ldots, p_k) \cdot \prod_{i=1}^{k} \varphi_{p_i} \right)$$

(15)

Where

$$\sum_{\sum_{i=1}^{k} p_i = n} a_k = \binom{n}{k}$$

Above for n=5, we have p(5)= 7 i.e. the partitions of 5: (5), (1, 4), (2, 3), (1, 1, 3), (1, 2, 2), (1, 1, 1, 2), and (1, 1, 1, 1, 1) with the weights $a_k$ given by (5), (5, 5), (5, 5), (5), (1) and the φ's given by Eq.(11). For n=6 we also have that

λ$_6$ = 6·φ$_6$ -(6·φ$_1$·φ$_5$ +6·φ$_2$·φ$_4$ +3·φ$_3$·φ$_3$ ) + (12·φ$_1$·φ$_2$·φ$_3$ + 6·φ$_1$·φ$_1$·φ$_4$ +2·φ$_2$·φ$_2$·φ$_2$ ) +
 - (6·φ$_1$·φ$_1$·φ$_1$·φ$_3$ + 9·φ$_1$·φ$_1$· φ$_2$·φ$_2$ ) + 6·φ$_1$·φ$_1$·φ$_1$·φ$_1$·φ$_2$ - 1·φ$_1$·φ$_1$·φ$_1$·φ$_1$·φ$_1$·φ$_1$ .

with the partition: (6), (1, 5), (2, 4), (3, 3), (1, 2, 3), (1, 1, 4), (2, 2, 2), (1,1 , 1, 3), (1,1,2,2), (1,1,1,1,2), (1, 1, 1, 1, 1, 1) (p(6)=11).

For the asymptotic and more about the partition function p(n) we refer to one of the fundamental work of Ramanujan [6, p.116].
In particular as

$$p(n) \sim \frac{e^{2\pi \cdot \sqrt{\frac{n}{6}}}}{(4 \cdot n \cdot \sqrt{3})} = \frac{e^{k \cdot \sqrt{n}}}{(4 \cdot n \cdot \sqrt{3})}$$

(16)

We have thus obtained an alternating binomial sequence -with "weighted cluster partitions" - whose values decrease rapidly to zero (the first two terms - the first as an upper bound - the first two as a lower bound- are expected to ensure the positivity of the λ's).
As a numerical experiment we first give on the Table below, the numerical values of the sequences furnishing the first 6 exact values of the λ's.

λ $_1$ = **0.023095708966** ..= 1+γ/2 – (½)·log(4π)
λ $_2$ = 0.0928791468860 – 0.000533411769594 =
    = **0.09234573511**..
λ $_3$ = 0.210844265482 - 0.00321766460904 +0.000012319522..=
    = **0.207638920396**.
λ $_4$ = 0.37949770799444 – 0.0108060650201 + 0.000099077079 - 2.84528115941·10$^{-7}$ =

= **0.368790470678**
λ $_5$ = 0.602359196025– 0.02727515800 80+0.000249373225556- 0.00000286056695512 +
　　+ 6.57137853964·10$^{-9}$ =
　　= **0.575330557247**
λ $_6$ = 0.8841219582034230 - 0.057948994032944 + 0.0014085173503118162574 +
　　− 0.00001554836720655352+ 7.928018702647952·10$^{-8}$–1.517706484293658·10$^{-10}$ =
　　= **0.8275659331**

(exact to some digits, the last one instead of 0.827566012.....[2]) .
We notice that the first term in each of the 6 sequences is an upper bound to the true value of the corresponding Li-Keiper coefficients and the first two terms is a lower bound to it as stated above.
As an example, for λ $_4$ , where the true value is λ $_4$= 0.368790479... [2], the upper bound is 0.37949..
and the lower one is 0.36869..
Thus: 0.36869... < λ $_4$ = 0.368790479...< 0.37949.
In Appendix 1, we illustrate the calculation for n=7, 8, 9, 10 where we found the inequality (the value of λ$_7$ is still taken from Ref [2]).
$$1.12078 < \lambda_7 = 1.12446011757 < 1.23044$$
and give also the expressions for the lower and the upper bound for n=8 , n=9 and n=10.
The general inequality for the Li-Keiper coefficients from this approach is formulated as:

$$n \cdot \varphi_n - \sum_{p=1}^{\left\lfloor \frac{n}{2} \right\rfloor} a_p \cdot \varphi_p \cdot \varphi_{n-p}$$

(17)

Where $\lfloor \ \rfloor$ is the floor function.
Since we are in the presence of an alternating decreasing sequence in absolute values, the successive contributions are positive confirming our inequality above for the upper and for
the lower bound along our lines.
In a numerical contest, an interesting numerical advanced experiment consists thus in using the first two "terms" above (the contribution of the "cluster" φ$_n$ (one block) and that of the product of only two "cluster contributions"(two blocks) as given in Eq. (17), where a$_p$ are the positive constants. The Table up to n=10 is given below (lower bound, true value, upper bound) with the plot.
The true values are taken from Maslanka [2].

| |
|---|
| λ $_1$ = **0.023095708966** ..= 1+γ/2 – (½)·log(4π) |
| λ $_2$ = 0.0928791468860 – 0.000533411769594 = <br> = **0.09234573522.** |
| 0.2076276009 < λ $_3$ = **0.20763892055** < 0.210844 26548 |
| 0.3686916430 < λ $_4$ = **0.36879047949** < 0.379497707994 |
| 0.5750840380 < λ $_5$ = **0.57542714461** < 0.602359196025 |
| 0.8261729642 < λ $_6$ = **0.8275660122** < 0.884121958203 |
| 1.1207865182 < λ $_7$ = **1.12446011757** < 1.230448804543 |
| 1.4573142832 < λ $_8$ = **1.46575567715** < 1.648279333807 |
| 1.8334010001 < λ $_9$ = **1.85091604838** < 2.145772946851 |
| 2.245753497 < λ $_{10}$ = **2.27933936319** < 2.732413708 |

　**Table 3**

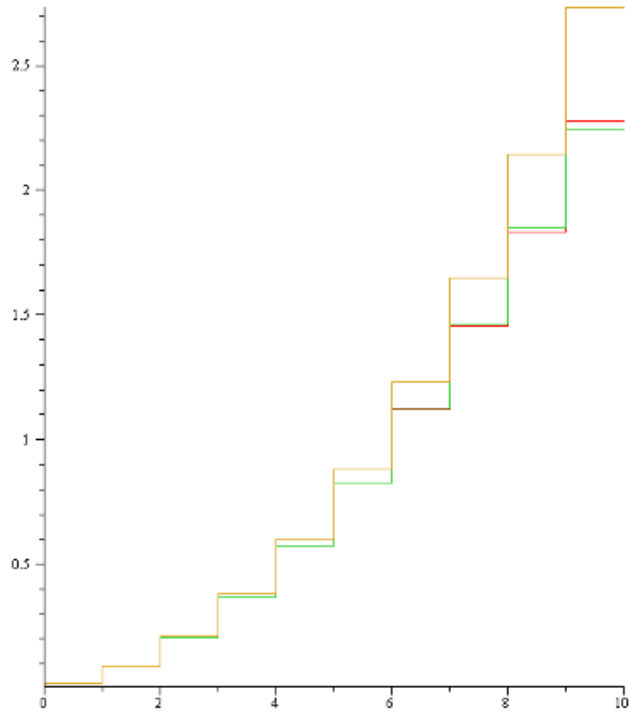

Fig.2 Plots of the lower, true values and upper bounds up to n=10 (From Table 3).

We have also pursued our numerical experiment and we have calculated our upper and lower bounds up to n=15.
In the next Fig. 3 we present the plot of our upper and lower bound from n=1 to n=15 together with a straight line i.e. 0.3·x + constant and the convex curve 0.1·x·log(x) + constant which illustrate the positivity property of our bounds given by Eq.(17). (If R.H. is true then the function will be given by 0.5· (x·log(x)) -c·x, for large integers x and c ~ 1.13 [2].

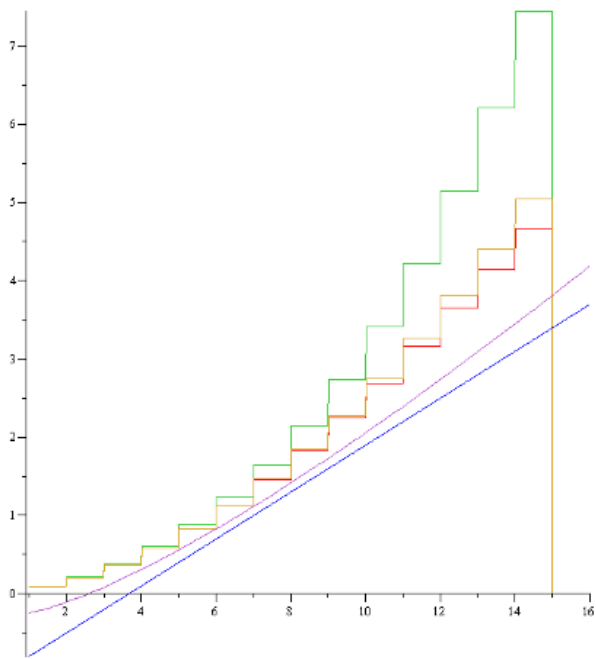

Fig.3 Plot of the lower bound, the true value, the upper bound and the function 0.1·x·log(x) + constant, the function 0.3·x + constant in the range (1..15) of x.

## 4. Shift of one unity in Xi(s) and the Glaisher-Kinkelin constant

The infinite set, given by Eq.(1) are in terms of the values of the successive derivatives of the Xi function at the right border of the critical strip, i.e. at s=1 (or at z=0, where z=1-1/s). Summation of the right hand sides of the set (which results in a strong cancellation at the left hand side) leads to a shift of one unity (from s=1 to s=2) in the argument of the Xi function and we obtain :

$$c = \sum_{n=1}^{\infty} \left(\frac{1}{n!}\right) \cdot \left(\frac{\xi^n(1)}{\xi(1)}\right) \cdot \frac{1}{\Gamma(n+1)} = \sum_{n=1}^{\infty} \left(\frac{1}{n!}\right) \cdot \left(\frac{\xi^n(1)}{\xi(1)}\right) \cdot \frac{(2-1)^n}{(n!)} =$$

$$= \left(\frac{\xi'(2)}{\xi(1)}\right) = \left(\frac{\pi}{3}\right) \cdot \left(\frac{\xi'(2)}{\xi(2)}\right)$$

(18)

Since $\xi(2) = (1/\pi) \cdot (\pi^2/6) = (\pi/3) \cdot (1/2) = (\pi/3) \cdot \xi(1)$, then:

$$c = \sum_{n=1}^{\infty} \left(\frac{1}{n!}\right) \cdot \left(\frac{\xi^n(1)}{\xi(1)}\right) \cdot \frac{1}{\Gamma(n+1)} = \left(\frac{\pi}{3}\right) \cdot \left(\frac{\xi'(2)}{\xi(2)}\right) =$$

$$= (\pi/3) \cdot ((1/2) \cdot (\gamma + \log(4 \cdot \pi) + 3) - 12 - \log(A))) = 0.072325988...$$

(19)

where A is the Glaisher-Kinkelin constant [7], given by
$\zeta'(s = -1) = 1/12 - \log(A) = -0.1654211937..$

(20)

Remark
Eq.(1) holds also for ξ replaced by $s.\pi^{-s/2} \cdot \Gamma(s/2)$ and the φ's related to the $\lambda_{trend}$'s as well as for the function $(s-1) \cdot \zeta(s)$ and the φ's related to the $\lambda_{tiny}$'s ($\lambda_{trend}(n) + \lambda_{tiny}(n) = \lambda(n) = \lambda_n$), where $\lambda_{tiny}(n)$ is the tiny part of the coefficient [2]. Then, the above constant split into the two contributions respectively given by:

$$(1/6) \cdot (-\gamma \cdot \pi - \pi \cdot \log(\pi) + \pi) = 0.4503357950 \text{ (trend)}$$

(21)

and

$$(1/6) \cdot ((12/\pi) \cdot \zeta'(s=2)) + 2 \cdot \pi) = -0.3780098064 \text{ (tiny)}$$

(22)

so that 0.4503357950 - 0.3780098064 = 0.0723259886 = c, the above constant.

Moreover, using the formula involving the number of zeros $N_o(n)$ [8], appearing in the binary representation of an integer n, we have - using the fundamental Formulas for γ and log(π) - [ 8 ], i.e.

$$\log\left(\frac{2}{\gamma}\right) = -\frac{1}{2} + \sum_{n=1}^{\infty} \left(\frac{N_1(n) + N_0(n)}{(2n) \cdot (2n+1) \cdot (2n+2)}\right)$$

(22)

and

$$\gamma = \frac{1}{2} + \sum_{n=1}^{\infty} \left( \frac{N_1(n) + N_0(n)}{(2n) \cdot (2n+1) \cdot (2n+2)} \right)$$

(23)

that

$$c = \left(\frac{\pi}{3}\right) \cdot \left( 1 + \frac{3}{2}\log(2) + 12 \cdot \zeta'(-1) + \sum_{n=1}^{\infty} \left( \frac{N_0(n)}{(2n) \cdot (2n+1) \cdot (2n+2)} \right) \right)$$

(24)

Numerically, using $N_o(n)$ up to n=32, we obtain:

$$c = \left(\frac{\pi}{3}\right) \cdot \left( 1 + \frac{3}{2}\log(2) + 12 \cdot \zeta'(-1) + \frac{1}{4 \cdot 5 \cdot 6} \ldots \right)$$

(25)

As noticed above, we note the strong cancellation in the left hand side and equal, now, to the linear function in the Li-Keiper coefficients which reads:

$$c = \left(\frac{\pi}{3}\right) \cdot \left( \frac{d}{dz}(\log(\xi(z))) \cdot \frac{d}{ds} \right)\bigg|_{z=\frac{1}{2}} = \left(\frac{\pi}{3}\right) \cdot \sum_{n=1}^{\infty} \lambda_n \cdot \left(\frac{1}{2}\right)^{n+1}$$

(26)

Here, using the first Li-Keiper coefficients [2] up to n=15 we obtain c = 0.072222733376... a value, as above, which may be compared with the exact value 0.072325...
An upper bound for c is also given by $\pi \cdot \lambda_1$ = 0.07255730.. [1] a value greater than the true value of only the amount $2 \cdot 10^{-4}$.

## 5. Concluding Remark and Riemann Hypothesis

In this work, we have presented an analysis of the infinite set of Equations for the Li-Keiper coefficients already obtained from their definition Eq.(1) [1], an alternative to the approach of a possible proof of the R.H using the Lee-Yang theorem on the zeros of the partition function for ferromagnetic Ising spin models with long range 2-body interactions in the presence of a magnetic field [1].
The set of Equations given by (1) is of course related to the Faà di Bruno's formula as well as to the complete Bell Polynomials for the φ's as well as to the partitions present in the structure of the algebraic expression for the Li-Keiper coefficients as functions of the φ's (the complete Bell polynomials $B_n(\lambda_1, \lambda_2, \ldots \lambda_n)$....) .
We are in presence of an alternating strong decaying sequence of contributions; in this context we have found a new formula for an upper and a lower positive bound to the Li-Keiper coefficients which involves only a single and a two cluster contributions.
The new formula has been checked up to n=15 and deserves future analysis both in a theoretical as well as in a numerical context. The constant c, related to the Glaisher-Kinkelin constant already obtained [1, 3], has been computed here using the formula for γ and for log(π) of Matiyasevich [8] and of Vacca and Sondow [9,10], in counting the zeros occurring in the binary representation of an integer n up to n=32, and also with the series with the λ's values up to $\lambda_{15}$ taken from Maslanka [2].
Finally, concerning the Riemann Hypothesis along ours lines we may subsume the following:

a. Dealing with Statistical Mechanics of spin-systems with two-body long range ferromagnetic interactions in the high temperature region - in connection with a truncation of the Xi function - we have found a possible proof of the R.H. by means of a lower "periodic" bound for the Li-Keiper coefficients given by the Riemann-Wave Background [1]. In this approach, we do not considered directly the behavior of nontrivial zeros and the main ingredient is the Lee-Yang theorem on the zeros of some of the partition functions of Statistical Mechanics. The lower bound is here "only" in the form of a "periodic" function of maximum value, which is equal to 4 [1, 3].

b. In the alternative approach along the above lines, we have used the definition of the Li-Keiper coefficients to obtain an infinite set of equilibrium equations for the Li-Keiper coefficients in connections with important functions like Bell polynomials, Faà di Bruno's formulas and partitions, i.e. a structure involving an alternating decaying sequence with the emergence of a new lower and upper bound for the Li-Keiper coefficients. To the best of our knowledge the bounds have not been considered before or along our lines.

These new bounds have been checked up to n=15 and they deserve of course further research.

The strong decaying alternating sequence related to the series expansion of the Xi function will be difficult to disprove. Moreover, the constant c, related to the Glaisher-Kinkelin [7] constant has been computed here in the binary system using important contributions of Matiyasevich [8].